\documentclass[12pt, a4paper]{amsart}

\usepackage{graphicx}
\usepackage{wrapfig}
\usepackage{comment}

\usepackage[
unicode=true,
plainpages = false,
pdfpagelabels,
bookmarks=true,
bookmarksnumbered=true,
bookmarksopen=true,
breaklinks=true,
backref=false,
colorlinks=true,
linkcolor = DarkBlue,
urlcolor  = DarkBlue,
citecolor = DarkRed,
anchorcolor = green,
hyperindex = true,
hyperfigures
]{hyperref}
\hypersetup{
pdftitle={Linear algebra, theory and algorithms},
pdfauthor={V.H. Mikaelian},
pdfsubject={Linear Algebra}
}

\usepackage[dvipsnames]{xcolor}
\colorlet{LightRubineRed}{RubineRed!70!}
\colorlet{Mycolor1}{green!10!orange!90!}
\definecolor{DarkRed}{HTML}{cc0000}
\definecolor{ChapterHeadColor}{HTML}{cc0000}
\definecolor{PartHeadColor}{HTML}{cc0000}
\definecolor{DarkBlue}{HTML}{0000cc}
\definecolor{QuoteColor}{HTML}{665665}

\usepackage[top=20mm, bottom=20mm, outer=26mm, inner=26mm, heightrounded
]{geometry}

\newcommand{\Z}{{\mathbb Z}}
\newcommand{\Q}{{\mathbb Q}}
\newcommand{\Co}{{\mathbb C}}

\usepackage[bitstream-charter]{mathdesign}

\DeclareSymbolFont{cmsymbols}{OMS}{cmsy}{m}{n}
\SetSymbolFont{cmsymbols}{bold}{OMS}{cmsy}{b}{n}
\DeclareSymbolFontAlphabet{\mathcal}{cmsymbols}

\theoremstyle{plain}
\newtheorem{Theorem}{Theorem}[section]
\newtheorem{Lemma}[Theorem]{Lemma}
  
\newtheorem{Corollary}[Theorem]{Corollary}  
  
\theoremstyle{definition}
  
\newtheorem{Example}[Theorem]{Example} 
\theoremstyle{remark}
\newtheorem{Remark}[Theorem]{Remark}

\numberwithin{equation}{section}

\usepackage[bitstream-charter]{mathdesign}

\DeclareSymbolFont{cmsymbols}{OMS}{cmsy}{m}{n}
\SetSymbolFont{cmsymbols}{bold}{OMS}{cmsy}{b}{n}
\DeclareSymbolFontAlphabet{\mathcal}{cmsymbols}

\begin{document}


\makeatletter
\@namedef{subjclassname@1991}{\textbf{MSC}}
\makeatother


\title[Embeddings using universal words]{Embeddings using universal words in the free group of rank $2$}

\author{V.\,H. Mikaelian
}

\begin{abstract}
For an arbitrary countable group given by its generators and defining relations  we discuss a specific method for its embedding into a certain 2-generator group. 
Our embedding explicitly lists the images of generators from the countable group in the 2-generator group, and from its relations it explicitly deduces some defining relations for the 2-generator group inheriting certain special features. 
The obtained method can be used in constructions of explicit embeddings of recursive groups into finitely presented groups. 
\end{abstract}

\maketitle

\setcounter{tocdepth}{3}

\let\oldtocsection=\tocsection
\let\oldtocsubsection=\tocsubsection
\let\oldtocsubsubsection=\tocsubsubsection
\renewcommand{\tocsection}[2]{\hspace{-12pt}\oldtocsection{#1}{#2}}
\renewcommand{\tocsubsection}[2]{\footnotesize \hspace{6pt} \oldtocsubsection{#1}{#2}}
\renewcommand{\tocsubsubsection}[2]{ \hspace{42pt}\oldtocsubsubsection{#1}{#2}}


\section{Introduction}
\noindent
The objective of this note is to suggest some simple rules for explicit embedding of any countable group $G$ given by its generators and defining relations into a $2$-generator group $T$ such that the defining relations of $T$ can easily be deduced from those of $G$, and they inherit certain features of the relations of $G$ needed for embeddings of recursive groups into finitely presented groups (see subsection~\ref{SU Preserving the structure}). 

By the well-known theorem of 
Higman, Neumann and 
Neumann an arbitrary countable group $G$ is embeddable into a $2$-generated group $T$ \cite{HigmanNeumannNeumann}. This result 
called in Robinson's textbook \cite{Robinson}
\textit{``probably the most famous of all embedding theorems''} was a starting step for
further research on embeddings into $2$-generator groups with related properties.
Typically, such research discusses cases when the embedding has an extra feature (is subnormal, is verbal, etc.), or when the group $T$ has a required property, including those inherited from $G$ (is soluble, is generalized soluble or generalized nilpotent, is linearly ordered, residualy has some property, is simple, etc.). For an outline of the topic see articles 
%
\cite{Dark}--\cite{On abelian subgroups}
and the literature cited therein.

In fact, the original embedding method of \cite{HigmanNeumannNeumann} and some of other embedding constructions cited above  already are explicit, and they do allow to compute the relations of $T$ based on the relations of $G$. However, we need a method that not only makes discovery of the relations of $T$ a simple, automated task, but also \textit{preserves certain features} in them required for study of embeddings of recursive groups into finitely presented groups (see references in \ref{SU Preserving the structure} below).

\medskip
We need the following notations to introduce the embedding. 
In the free group
$F_2=\langle
x,y
\rangle$ of rank $2$ consider some \textit{universal words}:
\begin{equation}
\label{EQ definition of a_i(x,y)}
a_i(x,y) = y^{(x y^i)^{\,2}\, x^{\!-1}} 
\!\! y^{-x} 
\!\!,\quad\quad i=1,2,\ldots
\end{equation}
(conventional notations $x^y=y^{-1}xy$,\; $x^{-y}=(x^{-1})^y$ are used here).
Assume a generic countable group $G$ is given by its generators and defining relations as: 
$$
G = \langle\, A \mathrel{|} R\, \rangle= \langle a_1, a_2,\ldots \mathrel{|} r_1, r_2,\ldots \,\rangle 
$$
where the $s$'th relation 
$r_s \in R$
is a word of length $k_s$ on  letters, say, 
$a_{i_{s,1}},\ldots,a_{i_{s,\,k_s}}\!\!\!\! \in A$.\,
If we replace in $r_s$
each $a_{i_{s,j}}$,\; $j=1,\ldots,k_s$,\; by the respective word 
$a_{i_{s,j}}(x,y)$ defined above, we get a new word
\begin{equation}
\label{EQ definition of r'_s}
r'_s (x,y)=
r_s\big(a_{i_{s,1}}\!(x,y),\ldots,a_{i_{s,\,k_s}}\!(x,y)\big)
\end{equation}
on just two letters $x,y$ in the free group $F_2$.
In these terms:

\begin{Theorem}
\label{TH universal embedding}
For any countable group $G = \langle a_1, a_2,\ldots \mathrel{|} r_1, r_2,\ldots \,\rangle $ the map $\gamma: a_i \to a_i(x,y)$,\; $i=1,2,\ldots$\,, defines an injective embedding of $G$ into the $2$-generator group 
$$
T_G=\big\langle x,y 
\;\mathrel{|}\;
r'_1 (x,y),\; r'_2 (x,y),\ldots\,
\big\rangle
$$
given by its relations 
$r'_s (x,y)$,\;
$s=1,2,\ldots$
\end{Theorem}

This is noting else but a new formulation of SQ-universality of $F_2$.
The proofs of Theorem~\ref{TH universal embedding}
and of its modification 
Theorem~\ref{TH universal embedding torsion-free} for torsion-free groups
occupy subsections~\ref{SU The universal generators}--\ref{SU Some simplification for torsion free groups} below.
Examples of applications with this embedding can be found in 
subsection~\ref{SU Examples of embeddings}. 

\medskip
We would like to stress the following case related to the question of Bridson and
de la Harpe mentioned as Problem 14.10 (b) in the Kourovka Notebook~\cite{kourovka}:
\textit{``Find an explicit embedding of $\Q$ in a finitely generated group; such a group
exists by Theorem IV in \cite{HigmanNeumannNeumann}''}.
The required explicit embedding of $\Q$ into a $2$-generator group $T$ was given in \cite{On a Problem on Explicit Embeddings of Q} in two manners, using free constructions and wreath products.
In the current note we add one more feature: 
$\Q$ can be explicitly embedded into such a group $T$ the defining relations of which
can also be explicitly listed. In Example~\ref{EX embedding of rational group}
we display an 
explicit embedding of $\Q$
into the $2$-generator group with directly given defining relations:
$$
T_\Q =\big\langle x,y 
\;\mathrel{|}\;
(y^s)^{(x y^s)^{\,2} x^{\!-1}} 
\!y^{-(x y^{s-1})^{\,2} x^{\!-1}}
\!\!,\;\;\;\; s=2,3\ldots
\big\rangle.
$$

Among other recent research on embeddings of $\Q$ 
into finitely generated groups we would like to briefly stress the following:
\cite{Darbinyan Mikaelian}
continues the \textit{linear order} relation of rational numbers in $\Q$ onto the whole $2$-generator group $T$, and it shows that the embedding can be \textit{verbal}.
\cite{On abelian subgroups}
mentions that $\Q$ can never be embedded into a finitely generated \textit{metabelian} group $T$ (see Section 7 in \cite{On abelian subgroups} and also \cite{Finiteness conditions for soluble groups}).
\cite{Adian Atabekyan}
provides an explicit verbal embedding of $\Q$ into a $2$-generator group $T=A_\Q(m,n)$ such that the \textit{center} of $T$ coincides with the image of $\Q$, i.e., $Z(T)\cong \Q$. 
One of the tasks in Problem 14.10 (a) \cite{kourovka} is to find an explicit embedding of $\Q$ into  a ``natural'' \textit{finitely presented} group.
\cite{The Higman operations and  embeddings} 
describes how Higman's procedure could be modified for a family of groups that includes $\Q$ to explicitly embed each of such groups into a $2$-generator finitely presented group.
And the first direct solution to the above problem appeared recently in \cite{Belk Hyde Matucci}. Moreover, one of the finitely presented groups constructed in \cite{Belk Hyde Matucci} is the group $T\!\mathcal{A}$ which is $2$-generator and simple.

\medskip
In the final subsection~\ref{SU Preserving the structure}
we refer to the main motivation that brought us to the study of embeddings in 
Theorem~\ref{TH universal embedding}
and 
Theorem~\ref{TH universal embedding torsion-free}:
the  constructive Higman embeddings \cite{Higman Subgroups in fP groups} of recursive groups into finitely presented groups.


\medskip
When the text of this note was uploaded to the arXiv.org I had an opportunity to discuss the topic with Prof. L.A. Bokut' who remarked interesting parallelism with 
\cite{Shirshov 58} where A.I. Shirshov 
constructed the elements 
$$d_k = \Big[a \circ \big\{[\cdots (a \underbrace{\circ \, b)\circ b \cdots ]\circ b}_k\big \}\Big]\circ (a\circ b),$$
\vskip-2mm
\noindent
$k=1,2,\ldots$, 
in the free associative algebra $A$ with two generators $a$ and $b$. Here $a\circ b$ denotes the Lie algebra product  $ab-ba$, and for details see \S 4 in  \cite{Shirshov 58}.
These elements freely generate a free Lie algebra $L(a,b)$ of countable rank.
They are used to define an embedding of any countably generated Lie algebra into a $2$-generator Lie algebra. The set 
$\{d_k \mathrel{|} k=1,2,\ldots \}$ is ``distinguished'' in the sense of \cite{Shirshov 56}. 
Meanwhile, the technique of our proof is very much different from 
\cite{Shirshov 58, Shirshov 56}, which makes this parallelism even more interesting.
See also \cite{Bokut 72} and  Remark~\ref{RE reference to HNN} below where we refer to a construction in the original work of Higman, Neumann and Neumann  \cite{HigmanNeumannNeumann}.

The current work is supported by the joint grant 20RF-152 of RFBR and SCS MESCS RA, and by the 18T-1A306 grant of SCS MESCS RA.


\section{References and some  auxiliary results}

\noindent
For general group-theoretical information we refer to \cite{Robinson, Kargapolov Merzljakov, Rotman}.
If $G = \langle\, A \mathrel{|} R \,\rangle$ is the presentation of the group $G$ by its generators $A$ and deifing relations $R$, then for an alphabet $B$ disjoint from $A$ and for any set $S\subseteq F_B$ of group words on $B$, we denote by
$\langle G, B \mathrel{|} S \,\rangle$ the group $\langle \,A \cup B \mathrel{|}  R \cup S \,\rangle$.
If $\varphi:G\to H$ is a homomorphism defined on a group $G=\langle g_1,g_2,\ldots \rangle$ by the images $\varphi(g_1)=h_1$,  $\varphi(g_2)=h_2,\ldots$\; of its generators, we may for briefness refer to $\varphi$ as the homomorphism \textit{sending}\, $g_1,g_2,\ldots$ \;to\; $h_1,h_2,\ldots$

\medskip

Our proofs in Section~\ref{SE Embeddings into 2-generator groups by ``universal'' generators} will be based on free constructions: the operations of free product of groups, of free product of a groups with amalgamated subgroup, and of HNN-extension of group by one or multiple stable letters (including the case with  infinitely many stable letters).
Background information on these constructions can be found in  \cite{Rotman, Lyndon Schupp, Bogopolski}.  Also, we refer to our recent notes \cite{A modified proof, Subvariety structures, The Higman operations and  embeddings} for specific notations that we also share here to write the 
free product
$G*_{\varphi} H = \langle G, H \mathrel{|} a=a^{\varphi} \text{ for all $a\in A$}\, \rangle$ of groups $G$ and $H$ with subgroups $A$ and $B$  amalgamated under the isomorphism $\varphi : A \to B$;\;
and the HNN-extension
$G*_{\varphi} t=\langle G, t \mathrel{|} a^t=a^{\varphi} \text{ for all $a\in A$}\, \rangle$
of the base group $G$ by the stable letter $t$ with respect to the isomorphism 
$\varphi : A \to B$ of the subgroups $A,B\le G$.\;
We also use HNN-extensions
$G *_{\varphi_1, \varphi_2, \ldots} (t_1, t_2, \ldots) = \langle G, t_1, t_2,\ldots \mathrel{|} a_1^{t_1}=a_1^{\varphi_1}\!\!,\; a_2^{t_2}=a_2^{\varphi_2}\!\!,\,\ldots\; \text{ for all $a_1\in A_1$, $a_2\in A_2,\ldots$}\, \rangle$
with multiple stable letters $t_1, t_2, \ldots$ \;
with respect to isomorphisms $\varphi_1: A_1 \to B_1,\; \varphi_2: A_2 \to B_2,\ldots$ for pairs of subgroups 
$A_1,B_1;\; A_2,B_2; \ldots$ in $G$.

\medskip
We are going to use certain subgroups of free products with amalgamated subgroups.  
Lemma~\ref{LE subgroups in amalgamated product} is a slight variation of Lemma 3.1 given on p.~465 of \cite{Higman Subgroups in fP groups} without a proof as \textit{``obvious from the normal form theorem for free products with an amalgamation''}. The proof can be found in subsection 2.5 of \cite{A modified proof}. 

\begin{Lemma}
\label{LE subgroups in amalgamated product}
Let $\Gamma = G *_\varphi \!H$ be the free product of the groups $G$ and $H$ with amalgamated subgroups
$A \le G$ and $B \le H$
with respect to the isomorphism 
$\varphi: A \to B$.
If $G', H'$ are subgroups of $G, H$ respectively, such that
for $A'=G'\cap A$ and $B'=H'\cap B$ 
we have  
$\varphi (A') = B'$, then for the subgroup $\Gamma'=\langle G',H'\rangle$  of $\Gamma$ and for the restriction $\varphi'$ of $\varphi$ on $A'$ we have:
\begin{enumerate}
\item $\Gamma' = G'*_{\varphi'} H'$,

\item $\Gamma' \cap A = A'$ and $\Gamma' \cap B=B'$,

\item $\Gamma' \cap G = G'$ and $\Gamma' \cap H = H'$.
\end{enumerate}
\end{Lemma}

If the amalgamated subgroups are trivial in a given free product with amalgamation, then that product is nothing but the ordinary free product of the same groups. Applying this observation to the groups $G'$ and $H'$ with trivial intersections $A'$ and $B'$ we get:

\begin{Corollary} 
\label{CO G*H free products}
In the notations of Lemma~\ref{LE subgroups in amalgamated product}:

\begin{enumerate}
\item if $A'=G'\cap A$ and $B'=H'\cap B$ both are trivial,
then $\Gamma'=\langle G',H'\rangle = G'*H'$,

\item if, moreover, $A'$ is a free group of rank $r_1$, and $B'$ is a free group of rank $r_2$, then $\Gamma'= F_r$ is a free group of rank $r=r_1+r_2$.
\end{enumerate}
\end{Corollary}

\section{Embeddings into $2$-generator groups by universal words}
\label{SE Embeddings into 2-generator groups by ``universal'' generators}

\subsection{The universal words in free group of rank $2$}
\label{SU The universal generators}

Let $F=F_A$ be a free group on a countable alphabet $A=\{a_1, a_2,\ldots\}$. Fix a new generator $a$, and in the free product $F'=F * \langle a \rangle$ 
set the isomorphisms of cyclic subgroups $\varphi_i:\langle a \rangle \to \langle a_i a \rangle$
sending $a$ to $a_i a$,\; $i=1,2,\ldots$
Define the respective HNN-extension:
$$
P\;=\;
F' *_{\varphi_1, \varphi_2,\ldots}
(t_1, t_2,\ldots)
\;=\;
\langle
F', t_1, t_2,\ldots \mathrel{|}\;
a^{t_i}= a_i a,\;\; i=1,2,\ldots\,
\rangle.
$$
Clearly, the stable letters $t_i$,\, $i=1,2,\ldots$\,,\,  generate in $P$ a free subgroup $X$ of countable rank.
It is simple to pick an auxiliary 2-generator group with a subgroup isomorphic to $X$: in the free group $Y=\langle
y,z
\rangle$ the elements $t_i'=y^i z^i$\!,\, $i=1,2,\ldots$\,, freely generate the subgroup $X'=\langle
t_1', t_2',\ldots
\rangle$ of countable rank. Amalgamating $X$ and $X'$
according to the isomorphism
$\psi$ sending $t_1,t_2,\ldots$ to $t'_1,t'_2,\ldots$
we get the group
$$
Q\;=\;
P *_{\psi} Y \;=\;
\langle\,
P,\; y,z \mathrel{|}
t_i=t_i',\;\; i=1,2,\ldots\,
\rangle.
$$
This group can already be generated by three elements $a,y,z$ because its generators
$$
a_i=a_i a \cdot a^{-1} = a^{t_i}a^{-1}
= a^{t_i'}a^{-1}
\quad \text{and}\quad\;\;
t_i=t_i'=y^i z^i
\!\!,\;\;\;\;\;
i=1,2,\ldots
$$
all are in $\langle a,y,z \rangle$.
By construction of $P$ no non-trivial power of 
$a$ is in $X$, and 
by construction of $Y$
no non-trivial power of 
$y$ is in $X'$, that is, the intersections  
$\langle a \rangle \cap X$ and  $\langle y \rangle\cap X'$  both are trivial. This by Corollary~\ref{CO G*H free products} implies that 
$\langle
a,y \rangle$ is a free subgroup of rank $1+1=2$ in $Q$. 
Introducing a new stable letter $x$ for the isomorphism $\pi:\langle  y,z\rangle \to \langle a,y\rangle$ sending $y,z$ to $a,y$, we construct:
\begin{equation}
\label{EQ formula of F_2}
F_2 =Q *_{\pi} x =
\langle
Q,x \mathrel{|}
y^x\!=a,\; z^x\!=y
\rangle
=\Big(\big((F* \langle a\rangle) *_{\varphi_1, \varphi_2,\ldots}
\!(t_1, t_2,\ldots)\big)*_{\psi} Y\Big)*_{\pi} x.
\end{equation}

This seems to cause confusion, as we earlier used $F_2$ to denote the free group of rank $2$ on the alphabet $\{x,y\}$. Let us verify that, in fact, \eqref{EQ formula of F_2} coincides with that free group. 

Firstly, $F_2$ is generated by $\{x,y\}$ because $a$ and $z$ can be expressed as some words on $x,y$:
$$
a=y^{x}=a(x,y),\quad
z=y^{x^{-1}}\!\!=z(x,y).
$$
Using these we can also 
express as words on $x,y$
\textit{all} the above discussed generators:
\vskip-4mm
$$
t_i = t_i' =\, y^i z^i
=\; y^i x y^i x^{-1}\!
=t'_i(x,y)=t_i(x,y),
$$
\vskip-5mm
\begin{equation}
\label{EQ formula for a_(x,y)}
a_i 
= a^{t_i}a^{-1} \!
= y^{x y^i \! x  y^i \! x^{\!-1}} 
\!\! y^{-x} 
= y^{(x y^i)^{\,2}\, x^{\!-1}} 
\!\! y^{-x}
=  \;
x  (y^{-i} x^{-1})^2   
y\,
(x\, y^i)^2    x^{-2} y^{-1} \! x =a_i(x,y).
\end{equation}
Secondly, $F_2$ is free  on $\{x,y\}$ because all the relations demanded by our ``nested'' free construction  (see the right-hand side in \eqref{EQ formula of F_2}) already hold on words
$
a(x,y),\, 
a_i(x,y),\, 
t_i(x,y),\,$
$ 
t'_i(x,y),\, 
y,\, 
z(x,y),\,  
x
$
in $F_2$, as it is easy to verify:
\begin{equation}
\label{EQ deducng everything from x, y}
\begin{split}
& \quad \;\;  a(x,y)^{t_i(x,y)}
= (y^x)^{\, y^i x y^i x^{-1}}\!\!\!
= y^{x y^i x y^i x^{-1}} \!\!\! y^{-x}\! y^{x}=
a_i(x,y)\, a(x,y),
\\
& t'_i(x,y)=t_i(x,y)
,\quad \quad
y^x = a(x,y),\quad\quad
z(x,y)^x = 
\big(y^{x^{-1}}\big)^x = y.
\end{split}
\end{equation}
No relations actually needed ``to bind'' $x$ with $y$, and so \eqref{EQ formula of F_2} is free on $x,y$. 

Clearly, $\bar F=\big\langle a_1(x,y), a_2(x,y), \ldots\big\rangle$ is an isomorphic copy of $F$ inside $F_2$,
and for any word $r(a_{i_1},\ldots,a_{i_k})\in F$ we have a word 
$$r' (x,y)=
r\,\big(a_{i_1}(x,y),\ldots,a_{i_k}(x,y)\big)\in \bar F \le F_2$$ obtained by replacing each $a_{i_j}$ by $a_{i_j}(x,y)$,\; $j=1,\ldots,k$.
%

\subsection{The embedding construction}
\label{SU The embedding construction}

Assume a  countable group $G$ is given as $G= \langle\, A \mathrel{|} R\, \rangle = \langle a_1, a_2,\ldots \mathrel{|} r_1, r_2,\ldots \,\rangle $ where the $s$'th relation 
$r_s(a_{i_{s,1}},\ldots,a_{i_{s,\,k_s}}\!)$
is a word on $k_s$ letters, as mentioned in  Introduction.
To such a relation we put into correspondence the word
$
r'_s (x,y)$
defined in $\bar F$ by 
\eqref{EQ definition of r'_s}. The set of such words for all $s=1,2,\ldots$ form a subset  $\bar R =\big\{r'_1 (x,y), r'_2 (x,y),\ldots\big\}$ with the normal closure $\bar N = \langle \bar R\rangle^{\bar F}$ in $\bar F$.
Since $\bar F \cong F$, we have $\bar N \cong N$ and $\bar F / \bar N \;\cong\; F/ N \;\cong\; G$.
Next define $\tilde N = \langle \bar R\rangle^{F_2}$ to be normal closure of $\bar R$ in the whole group $F_2$.

The natural homomorphism $\nu: F\to G \cong F/N$ is sending $a_i$ to $g_i=N a_i$, $i=1,2,\ldots$\,,\;
and we may consider the 
``updated'' isomorphisms 
$\varphi_i:\langle a \rangle \to \langle g_i a \rangle$ of cyclic subgroups in the free product 
$G * \langle a \rangle$
(for simplicity we do not introduce new letters for these $\varphi_i$).
In analogy with \eqref{EQ formula of F_2} we
build another ``nested'' free construction:
\begin{equation}
\label{EQ formula of H}
H=
\Big(\big((G * \langle a \rangle) *_{\varphi_1, \varphi_2,\ldots}
\!(t_1, t_2,\ldots)\big)*_{\psi} Y\Big)*_{\pi} x
\end{equation}
by using the new isomorphisms 
$\varphi_i$
and the same isomorphisms
$\psi$ and $\pi$
as used above to define $Q$ and $F_2$.

The natural homomorphism $\nu$ can be extended to a homomorphism $\bar \nu$ from the group $F_2$ onto $H$
by requiring $\bar\nu$ to agree with $\nu$ on $F$, and to fix each of the remaining generators
$a,t_1,t_2,\ldots,\,y,z,x$.

It turns out that the relations $\bar R$ already are enough to define the group $H$ because from \eqref{EQ formula of H} it is clear that all other equalities $a^{t_i} = g_i a$,\; $t_i=t'_i$, $y^x = a$,\; $z^x=y$ of $H$ do follow, like in \eqref{EQ deducng everything from x, y}, from representation of the generators via $x,y$.


Since $G$ trivially is embedded into $H$,  the subgroup $ (\bar F   \tilde N) / \tilde N $ of $F_2/\tilde N$ is isomorphic to $\bar F / \bar N\cong G$.
On the other hand we have $ (\bar F   \tilde N) / \tilde N \cong \bar F /(\bar F \cap \tilde N )$.
But since $\bar N$ is the kernel of the natural homomorphism from $\bar F$ to $\bar F /\bar N \cong G$, we get that $\bar F \cap \tilde N \le \bar N$.
Since also $\bar N \le \bar F$ and $\bar N \le \tilde F$, we have $
\bar F \cap \tilde N  = \bar N$. 

This construction together with~\ref{SU The universal generators} prove the following technical lemma:

\begin{Lemma}
\label{LE universal elements in F_2}
Let $F_2=\langle x,y\rangle$ be a free group of rank $2$ with  elements $a_i(x,y)$ defined in \eqref{EQ definition of a_i(x,y)} 
generating the subgroup $\bar F=\big\langle a_i(x,y) \mathrel{|} i=1,2,\ldots \big\rangle$ in $F_2$.
For any subgroup $N$ in $\bar F$ let $\bar N$ 
and $\tilde N$ denote
the normal closures of $N$ in $\bar F$ and in $F_2$ respectively. 
Then:
$$
\bar F \cap \tilde N  = \bar N.
$$
\end{Lemma}

Now we can conclude the proof of Theorem~\ref{TH universal embedding} as follows. 
Let the map 
$\gamma$,
the group 
$T_G=\big\langle x,y 
\mathrel{|}
r'_1 (x,y),\; r'_2 (x,y),\ldots\,
\big\rangle$,
and the relations
$r'_s (x,y)=
r_s\big(a_{i_{s,1}}\!(x,y),\ldots,a_{i_{s,k_s}}\!(x,y)\big)$
be those mentioned in the theorem.
Since the elements $a_1(x,y), a_2(x,y), \ldots$ generate an isomorphic copy $\bar F$ of $F$ in $F_2$, we have $G \cong \bar F / \bar N$.
Since $\bar N= \bar F \cap \tilde N$ by Lemma~\ref{LE universal elements in F_2},  then:
$$
G\;\cong\; \bar F / \bar N \;=\; \bar F / (\bar F \cap \tilde N)
\;\cong\;
(\bar F \tilde N) /\tilde N.
$$
But $\bar F \tilde N$ is in the whole $F_2$, and so $(\bar F \tilde N) /\tilde N$ clearly is a subgroup in $F_2 /\tilde N \;\cong\; T$.

\medskip
Theorem~\ref{TH universal embedding} provides a very easy way to embed a countable group $G = \langle a_1, a_2,\ldots \mathrel{|} r_1, r_2,\ldots\, \rangle $ into a $2$-generator group $T_G$ the relations of which are trivially obtained by just replacing 
in $r_1, r_2,\ldots$ all occurrences of the letters $a_1,a_2,\ldots$  by expressions $a_1(x,y),\; a_2(x,y),\ldots$

Notice that $T_G$ does depend on the particular presentation $G = \langle\, A \mathrel{|} R\, \rangle$, and for a different choice of $A$ and $R$ we may output another $2$-generator group. However, we do not want to note it as $T_{\langle\, A \,\mathrel{|} \,R\, \rangle}$ because this would bring to bulky notations in examples in  subsection 
\ref{SU Examples of embeddings}.

\subsection{Some simplification for torsion free groups}
\label{SU Some simplification for torsion free groups}

The isomorphisms $\varphi_i:\langle a \rangle \to \langle a_i a \rangle$
sending $a$ to $a_i a$ used in \ref{SU The universal generators} cannot, in general, be replaced by isomorphisms sending $a$ to $a_i$, $i=1,2,\ldots$,  because when $g_i \!\in\! G$ from \ref{SU The embedding construction} is an element of  \textit{finite} order, then $\langle a \rangle$ and $\langle g_i \rangle$ are \textit{not} isomorphic, and they can no longer be used as associated subgroups. This is the reason why we used $g_i a$ instead. But when 
$G$
is \textit{torsion-free}, this obstacle is dropped, and we can replace
$a_i a$ by $a_i$.
This allows to replace  $a_i(x,y)$ of \eqref{EQ formula for a_(x,y)} by a shorter word
\begin{equation}
\label{EQ definition of bar a_i(x,y)}
\bar a_i(x,y)
= a^{t_i}
=
y^{(x y^i)^{\,2} x^{\!-1}}
\!\!
,\quad\quad i=1,2,\ldots
\end{equation}
Replacing in $r_s$
each $a_{i_{s,j}}$ by  
$\bar a_{i_{s,j}}(x,y)$, we get other, shorter than $r'_s (x,y)$ word
$$
r''_s (x,y)=
r_s\big(\bar a_{i_{s,1}}\!(x,y),\ldots,\bar a_{i_{s,\,k_s}}\!(x,y)\big)
$$
on letters $x,y$ in the free group $F_2$. We have the following analog of 
Theorem~\ref{TH universal embedding}:

\begin{Theorem}
\label{TH universal embedding torsion-free}
For any torsion-free countable group $G\, = \,\langle a_1, a_2,\ldots \mathrel{|} r_1, r_2,\ldots \,\rangle $ the map $\gamma: a_i \to \bar a_i(x,y)$,\; $i=1,2,\ldots$\,, defines an injective embedding of $G$ into the $2$-generator group 
$$
T_G=\big\langle x,y 
\;\mathrel{|}\;
r''_1 (x,y),\; r''_2 (x,y),\ldots\,
\big\rangle
$$
given by its relations 
$r''_s (x,y)$,\;
$s=1,2,\ldots$
\end{Theorem}

Adaptation of the proof in \ref{SU The universal generators} and in \ref{SU The embedding construction} for this case is trivial.

\begin{Remark}
	\label{RE reference to HNN}
The reader may collate the above constructions with pages 252--254 in \cite{HigmanNeumannNeumann}. We used some ideas from there and from \cite{Higman Subgroups in fP groups}, but our proof is briefer, and we  produced shorter words $a_i(x,y)$. Compare them with  words 
$
e_i= 
a^{-1} b^{-1} a\, b^{-i} a\, b^{-1} a^{-1} b^{i}a^{-1} b\, a\, b^{-i} a\, b\, a^{-1} b^i
$
used in \cite{HigmanNeumannNeumann}.
And we have even shorter words $\bar a_i(x,y)$ for torsion-free groups.
\end{Remark}

\subsection{Examples of explicit embeddings}
\label{SU Examples of embeddings}

Here are some applications of the method  with
Theorem~\ref{TH universal embedding}
and with
Theorem~\ref{TH universal embedding torsion-free}.

\begin{Example}
\label{EX embedding of free abellian into 2-generator group}
The free abelian group $G=\Z^\infty$ of contable rank can be given as
$$
G = \big\langle a_1, a_2,\ldots \mathrel{|} [a_k,a_l],\; k,l=1,2\ldots \big\rangle
$$ 
by its relations $r_s=r_{k,l}=[a_k,a_l]$.
Since $G$ is torsion-free, we can use the shorter formula \eqref{EQ definition of bar a_i(x,y)} to 
map each $a_i$ to respective $\bar a_i (x,y)$. This defines an embedding of $\Z^\infty$ into the $2$-generator group:
$$
T_{\Z^\infty} =\big\langle x,y 
\;\mathrel{|}\;
\big[
y^{(x y^k)^{\,2} x^{\!-1}} 
\!\!\!\!\!,\,\,\,
y^{(x y^l)^{\,2} x^{\!-1}} 
\big] ,\;\;\; k,l=1,2\ldots
\big\rangle.
$$
\end{Example}

\begin{Example}
\label{EX embedding of rational group}
The additive group of rational numbers $G=\Q$ can be presented \cite{Johnson} as:
$$
G = \big\langle a_1, a_2,\ldots \mathrel{|} a_s^s=a_{s-1},\; s=2,3\ldots \big\rangle
$$ 
where the generator $a_i$ corresponds to the fractional number $\displaystyle {1 \over i!}$ with $i=2,3\ldots$
Rewrite each $a_s^s=a_{s-1}$ as 
$a_s^s\,a_{s-1}^{-1}$ and use the latter as the relation $r_s=r_s(a_{s-1},\, a_s)$ for each $s=2,3\ldots$
Since $G$ again is torsion-free, we can use the shorter formula \eqref{EQ definition of bar a_i(x,y)} to 
map each $a_i=\displaystyle {1 \over i!}$ to a $\bar a_i (x,y)$.
After easy simplification
$\bar a_i (x,y)=\big(y^{(x y^i)^{\,2} x^{\!-1}} \big)^i 
\big(y^{(x y^{i-1})^{\,2} x^{\!-1}}\big)^{-1}\!\!
=\,
(y^i)^{(x y^i)^{\,2} x^{\!-1}} 
y^{-(x y^{i-1})^{\,2} x^{\!-1}}
$
we get an embedding of $\Q$ into the $2$-generator group:
$$
T_\Q =\big\langle x,y 
\;\mathrel{|}\;
(y^s)^{(x y^s)^{\,2} x^{\!-1}} 
y^{-(x y^{s-1})^{\,2} x^{\!-1}}
\!\!\!\!,\;\; s=2,3\ldots
\big\rangle=\big\langle x,y 
\;\mathrel{|}\;
(y^s)^{x y^s x y} 
y^{-x y^{s-1}x}
\!\!,\;\; s=2,3\ldots
\big\rangle.
$$
\end{Example}

\begin{Example}
\label{EX embedding of Pruefer group}
The quasicyclic Pr\"ufer $p$-group $G=\Co_{p^\infty}$ can be presented as:
$$
G = \big\langle a_1, a_2,\ldots \mathrel{|}
a_1^p,\;\;\, a_{s+1}^p\!=a_s
,\;\; s=1,2\ldots \big\rangle
$$ 
where the generator $a_i$ corresponds to the 
primitive $(p^i)$'th root $\varepsilon_i$ of unity \cite{Kargapolov Merzljakov}.
As this group is \textit{not} torsion-free, we have to use the rather longer formula 
$a_i(x,y)$ from \eqref{EQ formula for a_(x,y)} as the image of $a_i$.
For the first relation 
$a_1^p$ of $G$ we get the new relation 
$a_1(x,y)^p=\big(y^{(x y)^{\,2}\, x^{\!-1}} 
\!\! y^{-x}\big)^p$.
Next, rewrite each $a_{s+1}^p=a_s$ as 
$a_{s+1}^p a_s^{-1}$\!\!,\, and use this as the relation $r_s=r_s(a_s,\, a_{s+1})$, with  $s=1,2\ldots$\;
The respective new relation will be 
$$
\big(y^{(x y^{s+1})^{\,2}\, x^{\!-1}} 
\!\! y^{-x}\big)^p 
\big(y^{(x y^s)^{\,2}\, x^{\!-1}} 
\!\! y^{-x}\big)^{-1}
=
\big(y^{(x y^{s+1})^{\,2}\, x^{\!-1}} 
\!\! y^{-x}\big)^p 
y^{x} y^{-(x y^s)^{\,2}\, x^{\!-1}}\!\!.
$$
And we have an embedding of $\Co_{p^\infty}$ into the $2$-generator group:
$$
T_{\Co_{p^\infty}} =\big\langle x,y 
\;\mathrel{|}\;\;
\big(y^{(x y)^{\,2}\, x^{\!-1}} 
\!\! y^{-x}\big)^p\!\!,\;\;\;\;
\big(y^{(x y^{s+1})^{\,2}\, x^{\!-1}} 
\!\! y^{-x}\big)^p 
y^{x} y^{-\,(x y^s)^{\,2}\, x^{\!-1}}
\!\!\!,\;\;\; s=1,2\ldots
\big\rangle.
$$
\end{Example}

\subsection{Usage in embeddings of recursive groups}
\label{SU Preserving the structure}

The main motivation why we needed the embeddings of 
Theorem~\ref{TH universal embedding}
and
Theorem~\ref{TH universal embedding torsion-free}
concerns study of constructive embeddings of recursive groups into finitely presented groups, i.e.,  \textit{constructive} Higman embeddings \cite{Higman Subgroups in fP groups}
(see, in particular, the references to embeddings of $\Q$ into finitely presented groups related to Problem 14.10 (a) \cite{kourovka} mentioned in Introduction).

One of the steps of the embedding for a recursive group 
$G = \langle a_1, a_2,\ldots \mathrel{|} r_1, r_2,\ldots \,\rangle$  (i.e., of a group with recursively enumerable relations $r_1, r_2,\ldots$) into a finitely presented group is the preliminary embedding to $G$ into a $2$-generator group $T=T_G$. And this group $T$ need also have 
recursively enumerable set of relations. 
The simple, automated embeddings that we built above do preserve that property.

Moreover, we need embeddings \textit{preserving special features} of relations. 
For the details we refer to \cite{The Higman operations and  embeddings}, and give just rough idea here. 
Each relation of a $2$-generator group $T$ can be coded by means of a certain sequence of integers. This allows to study the recursively enumerable sets of relations by means of certain sets of sequences of integers 
in \cite{Higman Subgroups in fP groups}. As we see in 
\cite{The Higman operations and  embeddings},
the embeddings of 
Theorem~\ref{TH universal embedding}
and
Theorem~\ref{TH universal embedding torsion-free} 
guarantee some close correlation between the 
relations of $G$ and those sets of sequences, which allows us to build constructive embedding of $G$ into a finitely presented group. 

In particular, compare Example~\ref{EX embedding of free abellian into 2-generator group} from this note to 
Example~3.1, 
Example~3.2 and 
Example 4.11 with ``abacus machine'' in 
\cite{The Higman operations and  embeddings}.


{\footnotesize
\vskip4mm
\begin{tabular}{l l}
Informatics and Appl.~Mathematics~Department
& 
College of Science and Engineering\\

Yerevan State University
&
American University of Armenia\\

Alex Manoogian 1 
& 
Marshal Baghramyan Avenue 40\\

Yerevan 0025, Armenia
&
Yerevan 0019, Armenia\\

E-mails:
\href{mailto:v.mikaelian@gmail.com}{v.mikaelian@gmail.com},
\href{mailto:vmikaelian@ysu.am}{vmikaelian@ysu.am}
$\vphantom{b^{b^{b^{b^b}}}}$

\end{tabular}
}

\end{document}